\documentclass[12pt]{article}
\usepackage{amsfonts}
\usepackage{amssymb,latexsym,amsmath, amscd}

\usepackage{amssymb,latexsym,amsmath}

\begin{document}

\def\no{\noindent}
\setlength{\parindent}{.25in}
\setlength{\textwidth}{6in}
\setlength{\oddsidemargin}{.25in}
\setlength{\evensidemargin}{.25in}
\setlength{\textheight}{9.5in}
\setlength{\headheight}{0in}
\setlength{\topmargin}{-.5in}

\newcommand{\qed}{{\unskip\nobreak\hfil
        \penalty50\hskip1em\hbox{}\nobreak\hfil
        $\square$\parfillskip=0pt\finalhyphendemerits=0 \par}\bigskip}

\newtheorem{dfn}{Definition}[section]
\newtheorem{add}[dfn]{Addendum}
\newtheorem{rem}[dfn]{Remark}
\newtheorem{thm}[dfn]{Theorem}
\newtheorem{mthm}[dfn]{Main Theorem}
\newtheorem{lem}[dfn]{Lemma}
\newtheorem{sublem}[dfn]{Sublemma}
\newtheorem{prop}[dfn]{Proposition}
\newtheorem{prob}[dfn]{Problem}
\newtheorem{ass}[dfn]{Assumption}
\newtheorem{classprob}[dfn]{Classical Problem}
\newtheorem{eigenvaluesofasum}[dfn]{Eigenvalues of a sum Problem}
\newtheorem{cor}[dfn]{Corollary}
\newtheorem{conj}[dfn]{Conjecture}
\newtheorem{ex}[dfn]{Example}
\newtheorem{ques}[dfn]{Question}
\newtheorem{techques}[dfn]{Technical Question}
\newtheorem{crit}[dfn]{Criterion}
\newtheorem{listof}[dfn]{List of Properties}
\newtheorem{conv}[dfn]{Convention}
\newtheorem{cons}[dfn]{Consequence}
\newtheorem{defn}[dfn]{Definition}
\newtheorem{fact}[dfn]{Fact}
\newtheorem{obs}[dfn]{Observation}
\newtheorem{warn}[dfn]{Warning}
\newtheorem{stabcrit}[dfn]{Stability Criterion}
\newenvironment{theorema}{\noindent{\bf Theorem A.}\em}{}
\newenvironment{theoremareal}{\noindent{\bf Theorem A (Real case).}\em}{}
\newenvironment{theoremb}{\noindent{\bf Theorem B.}\em}{}
\newenvironment{theoremcsymplectic}{\noindent{\bf Theorem C 
        (Symplectic case).}\em}{}
\newenvironment{theoremd}{\noindent{\bf Theorem D.}\em}{}

\let\lhd\vartriangleleft
\def\proof{\par\medskip\noindent{\it Proof. }}
\def\sketch{\par\medskip\noindent{\it Sketch of proof. }}
\def\lra{\longrightarrow}
\def\Lra{\Longrightarrow}
\def\ra{\rightarrow}
\def\Ra{\Rightarrow}
\def\CR{\curvearrowright}
\def\lh{\hookleftarrow}
\def\half{\frac{1}{2}}
\def\C{{\Bbb C}}
\def\R{{\Bbb R}}
\def\E{{\Bbb E}}
\def\H{{\Bbb H}}
\def\Z{{\Bbb Z}}
\def\K{{\mathcal K}}
\def\d{{\mathcal D}}
\def\P{{\Bbb P}}
\def\Q{{\Bbb Q}}

\def\B{{\Bbb B}}
\def\N{{\Bbb N}}
\def\F{{\mathcal F}}
\def\eps{\epsilon}
\def\al{\alpha}
\def\be{\beta}
\def\ga{\gamma}
\def\Ga{\Gamma}
\def\de{\delta}
\def\De{\Delta}
\def\Si{\Sigma}
\def\si{\sigma}
\def\L{{\mathcal L}}
\def\la{\lambda}
\def\La{\Lambda}
\def\Om{\Omega}
\def\om{\omega}
\def\D{\partial}
\def\hook{\hookrightarrow}
\def\embed{\hookrightarrow}
\def\8{\infty}
\def\<{\langle}
\def\>{\rangle}
\def\e{\sim}
\def\BE{\begin{equation}}
\def\EE{\end{equation}}
\def\LL{{\mathcal L}}
\def\geo{\partial_{\infty}}
\def\tits{\partial_{Tits}}
\def\tangle{\angle_{Tits}}
\def\ov{\overrightarrow}
\def\ol{\overline}
\def\grad{\mathop{\hbox{grad}}}
\def\8{\infty}
\newcommand{\restr}{\mbox{\Large \(|\)\normalsize}}
\def\oo{{\cal O}}
\def\too{\tilde{\cal O}}
\def\3{\ss}
\def\Aql{A_{\mathfrak{q}}(\lambda)}

\def \pol{{\mathcal P}_n(\mathfrak g)}

\newcommand{\bbCPm}{\bbC\mathbb{P}^m}

\def\check{\centerline{\bf CHECK THIS!}}
\def\gap{\centerline{\bf GAP!}}

\def\goth{\mathfrak}

\overfullrule=0pt

\title{The cohomology with local
coefficients of compact hyperbolic manifolds - long version}
\author{John Millson 
\footnote{Partially supported by NSF grant number DMS 0104006}
 }
\date{\small{\it{ Dedicated to M.S Raghunathan }}}
\maketitle

\begin{abstract}
We extend the techniques developed by Millson and Raghunathan in \cite{MillsonRaghunathan}
to prove nonvanishing results for the cohomology of compact
arithmetic quotients $M$ of hyperbolic $n$-space $\mathbb{H}^n$
with values in the local coefficient systems associated to  finite
dimensional irreducible representations of the group $SO(n,1)$. We prove that
all possible nonvanishing results compatible with the vanishing
theorems of \cite{VoganZuckerman} can be realized by sufficiently deep
congruence subgroups of the standard cocompact arithmetic examples.

\end{abstract}

\section{Introduction}
The purpose of this paper is to extend the techniques of \cite{MillsonRaghunathan}
to prove nonvanishing results for the cohomology of certain 
arithmetic quotients $M$ of hyperbolic $n$-space $\mathbb{H}^n$
with values in the local coefficient systems associated to a finite
dimensional representations $W$ of the group $SO(n,1)$. In the compact case
we replace the technique of intersecting pairs of totally-geodesic
submanifolds considered there by intersecting the same submanifolds
but now each is  equipped with a {\em local coefficient}, i.e. a parallel
section of the restriction of the associated local coefficient system $\widetilde{W}$ 
restricted to the cycle (or equivalently a nonzero vector in $W$
fixed under the fundamental group of the submanifold). In \S 3.3.1 
we explain how such cycles correspond to an especially simple class of
cycles in the Eilenberg MacLane complex $C_.(\Gamma) \otimes W$.
To define an intersection pairing of such cycles one needs a pairing on the 
local coefficients. The key fact that it can be arranged that the two 
manifolds intersect in 
a single connected component we borrow from the earlier papers
\cite{MillsonRaghunathan}, \cite{JohnsonMillson} and \cite
{FarrellOtanedaRaghunathan}. The 
remaining problems are to find for which $W$ the required local coefficients
exist and then to verify that the coefficient pairings applied to the local
coefficients are nonzero.

In order to state our main theorem we need a definition.

\begin{dfn}
Let $\mu = (b_1,b_2,\cdots,b_m)$ be a dominant weight for $SO(n,1)$
where $m = [\frac{n+1}{2}]$ and we use the usual Cartesian coordinates
on the dual of the Cartan (so the positive roots are the sums and differences
of the coordinate functionals in case $n+1$ is even and the
coordinate functionals are the short positive roots if $n+1$ is odd).
Then we define $i(\mu)$ to be the number of nonzero entries in $\mu$.
\end{dfn}

In what follows we will let $\Gamma$ denote a member of
the special family of cocompact torsion-free lattices of $SO(n,1)$
which are  congruence subgroups of some level $\mathfrak{b}$
of the group of units of the form $f(x)$ in $n+1$ real variables given by
$$f(x) = x_1^2 + \cdots + x_n^2 - \sqrt{m} \ x_{n+1}^2$$ 
for some square-free positive integer
$m$ and a sufficiently  large ideal $\mathfrak{b}$ in the integers of $\Q(\sqrt{m})$. 
We will let  $\Phi$ denote a general cocompact torsion-free lattice of $SO(n,1)$.

\begin{thm}\label{bigtheorem}
Let $W$ be the irreducible representation of $SO(n,1)$ with highest
weight $\mu$. Then we have
\begin{enumerate}
\item If $n = 2m-1$  and all entries of $\mu$ are nonzero 
then  $H^p(\Phi,W) = 0$ for all cocompact lattices
$\Phi$ of $SO(n,1)$ and all $p$.
\item
For all other $W$ we have 
$$p\notin \{i(\mu), i(\mu) +1, \cdots, n- i(\mu) \} \Rightarrow
H^p(\Phi,W) = 0.$$

\item For all other $W$ we have (for $\mathfrak{b}$ sufficiently large 
depending on $W$)
$$p\in \{i(\mu), i(\mu) +1, \cdots, n- i(\mu) \} \Rightarrow
H^p(\Gamma,W) \neq 0.$$
\end{enumerate}
\end{thm}

\begin{rem} The vanishing results are (essentially) due to
Vogan and 
\newline
Zuckerman,\cite{VoganZuckerman}. A careful statement (without computational
details)
of what \cite{VoganZuckerman} implies for $SO(n,1)$ can be found
in \cite{RohlfsSpeh}, \S 1.3. 
Roughly half the above nonvanishing
results for the cocompact case were proved by J.- S. Li in \cite{Li}.
The rest are new as is the formulation of the theorem.
See the end of the Introduction for a more complete discussion
including a discussion of related results in the noncocompact case.
If we consider the more general case of lattices $\Phi$ in $Spin(n,1)$ and 
representations of $Spin(n,1)$ the analogue of the above statements 1. and 2.
above still hold by \cite{VoganZuckerman}. As for statement 3. one necessarily 
has $i(\mu) = m$ for otherwise one does not have a genuine representation of 
$Spin(n,1)$ and the only case of interest is $n = 2m$ .
Let $\Phi$ be a torsion-free cocompact lattice in  $Spin(2m,1)$ and $W$
have highest weight $\mu$ satisfying $i(\mu) =m$.  Then one obtains 
$H^m(\Phi,W) \neq 0$ by the usual Euler characteristic argument.
\end{rem}

We can also prove nonvanishing results for cup-products of cohomology groups
with local coefficients. We give only one such example here. 

The cohomology algebra with trivial coefficients $H^*(\Gamma,\R)$ acts on the
cohomology groups with local coefficients. We remind the reader that
it was proved in \cite{MillsonRaghunathan} that for the above cocompact lattices 
$$H^p(\Gamma,\R) \neq 0, 0 \leq p \leq n.$$

In the last section of this paper we will prove a general result which
implies

\begin{thm}
Suppose $W$ is an irreducible finite-dimensional representation of $SO(n,1)$
and $H^q(\Gamma,W) \neq 0$. Suppose $p \geq 1$ and that $p+q\leq [\frac{n}{2}]$.
Then the cup-product
$$H^p(\Gamma,\R) \otimes H^q(\Gamma,W) \to H^{p+q}(\Gamma, W)$$
is a nonzero map.
\end{thm}

\begin{rem}
We can use the previous theorem to explain why the set of positive integers
$k$ such that $H^k(\Gamma,W) \neq 0$ is
an unbroken string of integers. Indeed by operating on $H^*(\Gamma,W)$ by
the cohomology with trivial coefficients and applying the previous theorem
we have $H^{i(\mu)}(\Gamma,W) \neq 0 \Longrightarrow 
H^k(\Gamma,W) \neq 0,i(\mu) \leq k \leq [\frac{n}{2}]$. 
But then by Poincar\'e duality we have $H^k(\Gamma,W) \neq 0,
i(\mu) \leq k \leq n - i(\mu)$.

\end{rem}

It is important to observe that there is considerable overlap between
this theorem and the work of Jian-Shu Li, \cite{Li}. In 
\cite{Li}, Professor
Li proves  nonvanishing results for  cohomology groups with
local coefficients for cocompact lattices in the classical  
groups $SO(p,q), Sp(p,q), SU(p,q), Sp_{2m}(\mathbb{R})$ and $SO^*(2n)$.
When specialized to the case of $SO(n,1)$ his results give roughly
half of the nonvanishing results in our Theorem \ref{bigtheorem}
in that he proves nonvanishing for those $p$ in the range
$i(\mu) \leq p < \frac{n-1}{4}$ (and the Poincar\"e dual dimensions).
We also note that the special case  $H^1(\Phi,W) \neq 0 \Leftrightarrow 
i(\mu) =1$ of the above theorem is a consquence of the results of
\cite{Raghunathan} and \cite{Millson1985}. The result that
$i(\mu) > 1 \Rightarrow H^1(\Phi, W)=0$ for all $\Phi$ is a special case
of Raghunathan's Vanishing Theorem, \cite{Raghunathan} and the
result that $i(\mu)=1 \Rightarrow H^1(\Gamma, W) \neq 0$ 
was proved in \cite{Millson1985}.

We also need to mention the large amount of work on the noncompact case.
In this case there is a large supply of cohomology coming from the Borel-
Serre boundary, \cite{BorelSerre}. The idea of using Eisenstein series to promote
boundary classes to classes inside began with the fundamental papers
of Harder \cite{Harder1973} and \cite{Harder1975} and nonvanishing
theorems for this case analogous to the one we prove for the compact case
were proved by  Harder, loc. cit., Rohlfs and Speh, \cite{RohlfsSpeh} and Speh,
\cite{Speh}. For the
noncompact case one has the harder problem of finding cuspidal
classes. Progess on this for the case treated in this paper is to
be found in \cite{RohlfsSpeh}. More general results
that have some application to the case considered here can be found in
\cite{BorelLabesseSchwermer},\cite{LiSchwermer} and \cite{Clozel}. 
It is possible to use the intersection-theoretic methods of this paper
to prove nonvanishing results for the finite-volume noncompact case.
We will discuss this in later work.

In joint work with Jens Funke, \cite{FunkeMillson2}, we will show that the results of
\cite{KudlaMillson} concerning the connection between the theta
lifting and cycles in cocompact quotients of the symmetric spaces
of $SO(p,q)$ and $SU(p,q)$ generalize to local coefficients. In particular
the Poincar\'e dual harmonic forms of some of the cycles 
with local coefficients considered here can be 
obtained from an appropriate theta-lifting of Siegel modular forms of 
genus equal to the codimension of the cycles.

I would like to thank Jeffrey Adams for a number of helpful conversations
about the material in this paper, especially about the section concerned
with applying the Vogan-Zuckerman Theorem. I would also like to thank
Armand Borel for pointing out an error in an earlier version of this
paper. I would like to thank
the second referee for reminding me about the work of Jian-Shu Li, \cite{Li}, 
and Jian-Shu Li for explaining the relation between his results
and our Theorem \ref{bigtheorem}. Finally, I would like to thank the
first referee for pointing out the reference \cite{RohlfsSpeh} which drew
my attention to the noncompact finite volume case and for his suggestions
for the format of the paper.

This paper is dedicated to M.S.Raghunathan on the occasion of his sixtieth
birthday.  It is different from the one that I presented at the conference
at the Tata Institute in his honor in December of $2001$. It was felt that 
paper was too long
for the conference volume and it appeared difficult to split it up
in a satisfactory fashion. However this paper is perhaps more appropriate
in that it depends in an essential way on the techniques Raghunathan
and I developed in \cite{MillsonRaghunathan}. It is my pleasure to 
acknowledge the great impact that Raghunathan had at the beginning of my 
career, by providing a critical insight in \cite{Millson1976}
and by the collaboration \cite{MillsonRaghunathan}. 

\section{The vanishing theorem}

The goal of this section is to prove the following theorem. 
Let $m= [\frac{n+1}{2}]$. Let $W$ be the irreducible representation of 
$Spin(n,1)$ with highest weight $\mu= (b_1,b_2,\cdots,b_m)$. Recall
that $i(\mu)$ is the number of nonzero entries of $\mu$.

\begin{thm} \label{vanishingtheorem}
Let $\Phi$ be a cocompact lattice in $SO(n,1)$. Then
\begin{enumerate}
\item If $n$ is odd and $i(\mu) = l$ then
$$H^p(\Phi,W) =0, \  \text{for all p}.$$
 
\item In all other cases
$$H^p(\Phi,W) =0, p\leq i(\mu) - 1. $$
\end{enumerate}
\end{thm}

The Theorem \ref{vanishingtheorem} follows by working out what the results of
Vogan and Zuckerman in \cite{VoganZuckerman} imply for the case in hand.
We will establish some notation and then review the theory of
Vogan and Zuckerman. 

First we will use the following convention to deal with the problems of
denoting a Lie group or Lie algebra and its complexification. We will use Roman
letters e.g. $G,K,L$ to denote real Lie groups. We will use the corresponding 
gothic letters with a subscript $0$ e.g $\mathfrak{g}_0, \mathfrak{k}_0,
\mathfrak{l}_0$ to denote their Lie algebras and the same symbol without
the subscript $0$ to denote the complexification of the corresponding
real Lie algebra.
Furthermore, a gothic letter without a subscript $0$ e.g. $\mathfrak{q}$
or $\mathfrak{u}$ will denote a {\em complex} Lie algebra.

We let $V$ be an $n+1$-dimensional real vector space equipped with a
symmetric bilinear form of signature $(n,1)$ . We choose a basis
$\{e_1,\cdots,e_{n+1}\}$ such that
$(e_i,e_j) = 0, i \neq j$,  $(e_i,e_i) =1, 1\leq i \leq n$ and
$(e_{n+1},e_{n+1}) = -1$. We will need to use a Witt basis for
the complexification $V\otimes \C$.

In case $n+1$ is odd put $n =2m$. Define
$u_i = e_i - \imath e_{m+i}$ and $v_i =  e_i + \imath e_{m+i},
1\leq i \leq m$.
We arrange the basis in the order $\{u_1,\cdots,u_m,e_{n+1},v_m,\cdots,v_1\}$.

In case $n+1$ is even, put $n+1=2m$. Define $u_i$ and $v_i$ for 
$1\leq i \leq m-1$ as before. Define $u_m$ and $v_m$ by
$u_m = e_m - e_{2m}$ and $v_m = e_m + e_{2m}$.
We then define our ordered basis to be $\{u_1,\cdots,u_m,v_m,\cdots,v_1\}$.

In both cases we will call the above basis a Witt basis (even though
$(u_i,v_i) = 2, 1\leq i \leq m$).

We define a Cartan involution $\theta:V \to V$ by
$$\theta(e_i) = \begin{cases}  e_i, & 1 \leq i \leq n \\
-e_{i}, & i = n+1.
\end{cases}.$$
The centralizer of $\theta$, to be denoted $K$, is then a maximal compact
subgroup isomorphic to $O(n)$. We will let $K^0$ denote the connected
component of the identity of $K$. We let $o$ be the point in $\mathbb{H}^n$
fixed by $K$.

We let $\mathfrak{g}_0$ denote the Lie algebra of $G=SO(n,1)$
and $\mathfrak{k}_0$ denote the Lie algebra of $K^0$. We let
$\mathfrak{p}_0$ be the orthogonal complement to $\mathfrak{k}_0$
for the Killing form. Thus we have a canonical isomorphism
$$\mathfrak{p}_0\cong T_o(\mathbb{H}^n).$$

We represent $\mathfrak{g}$ as $n+1$ by $n+1$ matrices using the Witt basis.
We let $\mathfrak{h}$ be the Cartan subalgebra of  $\mathfrak{g}$
consisting of diagonal $n+1$ by $n+1$ matrices with diagonal entries
$(x_1, \cdots,x_l,-x_l,\cdots,-x_1)$ in case $n+1$ is even and
$(x_1, \cdots,x_l,0,-x_l,\cdots,-x_1)$ in case $n+1$ is odd.

In what follows it will be convenient to use the canonical isomorphism
$\phi:\bigwedge^2V \to \mathfrak{so}(V)$ for $V$ a (real or complex) vector
space equipped with a nondegenerate symmetric bilinear form $(\ ,\ )$.
The isomorphism is determined by the formula
$$\phi(x \wedge y) (z) = (x,z)y - (y,z)x.$$
We will use $\phi$ to identify $\bigwedge^2V$ and $\mathfrak{so}(V)$
henceforth.

We have 
 $$\mathfrak{h} = \oplus_i^m \C  \ u_i\wedge v_i.$$
In the case that $n+1$ is odd we have $\mathfrak{h} \subset \mathfrak{k}$.
In the case that $n+1$ is even we have 
$$\mathfrak{h} \cap \mathfrak{p} = 
\C \  u_m \wedge v_m = \C \  e_m \wedge e_{2m}= \mathfrak{a}.$$
Here $\mathfrak{a}$ is the complexification of the Cartan subspace
$\mathfrak{a}_0 = \R  \ e_m \wedge e_{2m}$ of $\mathfrak{p}_0$. 
We set $\mathfrak{t} := \mathfrak{h} \cap \mathfrak{k}$ whence
(if $n+1$ is even)

$$ \mathfrak{h} = \mathfrak{t} \oplus \mathfrak{a}.$$   
  
\subsection{The representations $\Aql$.}

In \cite{VoganZuckerman} the authors constructed all unitary representations
$\pi$ of a semisimple Lie group $G$ with nonzero continuous cohomology in a 
finite dimensional representation $W$ of $G$. These representations
are parametrized by pairs $(\mathfrak{q}, \lambda)$ consisting of a 
$\theta$--stable parabolic subalgebra $\mathfrak{q}=\mathfrak{l}
\oplus \mathfrak{u}$ and a unitary
one-dimension representation $\lambda$ of $L$. 
They will be denoted $\Aql$. We will also use $\lambda$
to denote the corresponding representation  of $\mathfrak{l}_0$. 
In what follows we let $L^c$
be the compact form of the complexification $L_{\C}$ of $L$.

Let $R = dim (\mathfrak{u} \cap \mathfrak{p})$ and
$S = dim (\mathfrak{u} \cap \mathfrak{k})$. The following fundamental
lemma is proved in \cite{BorelWallach}, the first statement is Theorem 5.2,
pg. 222 and the second (Wigner's Lemma) is Corollary 4.2, pg.26. 

\begin{thm}
In order that there exist a cocompact lattice $\Gamma \subset G$ such that
$H^p(\Gamma,W)\neq 0$ it is necessary that there exists
an irreducible unitary representation $\pi$ such 
that the 
$$H^p_{cont}(G,\pi \otimes W^*) \neq 0.$$
Furthermore the nonvanishing of $H^p_{cont}(G,\pi \otimes W^*) \neq 0$
implies that the infinitesimal characters of $\pi$ and $W$ coincide.
\end{thm}

To obtain Theorem \ref{vanishingtheorem} we have only to examine which
representations $\Aql$ have the same infinitesimal character as that
of $W$ and calculate the lowest dimensions in which they have continuuous
cohomology. 

We have 
\begin{lem}\label{VZ}
\begin{enumerate}
\item The infinitesimal character of $\Aql$ is $\lambda +\rho(\mathfrak{g})$
\item Let $W$ be the finite dimensional representation with
infinitesimal character $\lambda +\rho(\mathfrak{g})$ (so $W$ has highest
weight $\lambda$).
Then the representation $A_{\mathfrak{q}}(\lambda) \otimes W^*$ has cohomology  degrees $R +j$ 
 for those $j$ such that $H^j(L^c/L^c \cap K,\R) \neq 0$.
\item In particular the only finite dimensional representation $W$ for which
$H^k(A_{\mathfrak{q}}(\lambda)\otimes W^*)\neq 0$ for any $k$ is the representation $W_{\lambda}$
with highest weight $\lambda$ and the
lowest degree $k$ in which $A_{\mathfrak{q}}(\lambda) 
\otimes W_{\lambda}^*$ can have  nonzero cohomology
is in degree $R$.
\end{enumerate}
\end{lem}

\proof
The first statement is (b) of Theorem 5.3 of \cite{VoganZuckerman}
and the second in contained in Theorem 5.5 of \cite{VoganZuckerman}.
\qed

Let $\Delta \in \mathfrak{h}$ be the dual cone to the cone of positive roots and
$x \in \Delta$. We recall the definition of  the standard
parabolic subalgebra $\mathfrak{q}(x)$ associated to $x$.
We have

\begin{align*}
 \mathfrak{q}(x) = \  &\mathfrak{l}(x) \oplus \mathfrak{u}(x)\\
 \text{with}\\
 \mathfrak{l}(x) = \ & \mathfrak{h} \oplus \sum\limits_{\<\alpha,x\>=0}
 \mathfrak{g}_{\alpha}
\  \text{and} \ 
 \mathfrak{u}(x) =  \sum\limits_{\<\alpha,x\>>0}
\mathfrak{g}_{\alpha}
\end{align*}

We also define the standard parabolic $\mathfrak{q}(\lambda)$ associated to
$\lambda$ in the weight cone (so $\lambda \in \mathfrak{h}^*$ and
$\<\lambda ,\alpha^{\vee}\> \geq 0$ for all positive coroots $\alpha^{\vee}$).
The algebra is defined in the same way as $\mathfrak{q}(x)$ be replacing $x$
by $\lambda$ and $\<\ ,\ \>$ by the (dual) Killing form $(\ ,\ )$.

In what follows we will let $\Delta(x)$, resp. $\Delta(\lambda)$, denote the subset
of the {\em positive} roots that are orthogonal to $x$, resp. $\lambda$.

\begin{dfn} We will say a character $\lambda \in \mathfrak{h}^*$ is compatible
with a parabolic subalgebra $\mathfrak{q}$ if $\lambda$ extends to
a character of $\mathfrak{l}$ or equivalently if $\lambda$ annihilates
the commutator $[\mathfrak{l},\mathfrak{l}]$.
\end{dfn}

Note that the definition of the representation $\Aql$ assumes that
$\mathfrak{q}$ and $\lambda$ are compatible.

We have 

\begin{lem}\label{reducetostandard}
If $\lambda$ is compatible with a standard parabolic subalgebra 
$\mathfrak{q}(x)$ then
$$ R(\mathfrak{q}(x)) = dim(\mathfrak{u}(x)\cap \mathfrak{p}) \geq 
dim(\mathfrak{u}(\lambda) \cap \mathfrak{p})= 
R(\mathfrak{q}(\lambda)).$$ 
\end{lem}

\proof
We claim that $\Delta(x) \subset \Delta(\lambda)$.
Indeed suppose that $\alpha \in \Delta(x)$. Then
$h_{\alpha} = [x_{\alpha},y_{\alpha}]$ with $x_{\alpha},y_{\alpha}  \in 
\mathfrak{l}(x)$. Hence $h_{\alpha} \in [\mathfrak{l}(x),
\mathfrak{l}(x)]$. Since by assumption $\lambda$ annihilates
this commutator subalgebra we have $\lambda(h_{\alpha})=0 \Rightarrow 
(\lambda,\alpha) =0 \Rightarrow \alpha \in \Delta(\lambda)$. The claim
follows.
Thus  $\Delta(\lambda)^c \subset \Delta(x)^c$.
Equivalently, for a positive root $\beta$ 
we have $(\beta,\lambda) > 0
\Rightarrow \<\beta,x\> > 0$ and
consequently
$$\mathfrak{u}(\lambda) \subset \mathfrak{u}(x)$$
and the lemma follows.
\qed

We will first deal with the case of $SO(2m,1)$.

\subsection{The case of $SO(2m,1)$}
We have $l=m$. 
Now suppose we are given $W$ with highest weight
$(b_1,\cdots,b_m)$. We assume $i(\lambda) = r$ and we set $s = m -r$.

Theorem \ref{vanishingtheorem} for the case $n=2m$ will follow from

\begin{lem}\label{Rlemma}
Suppose that $\lambda$ is as above. Then we have
\begin{enumerate} 
\item $L(\lambda) =C \times SO(2s,1)$ where $C \subset K$.
\item The set $\{ u_i \wedge e_{2m+1}: 1\leq i \leq r \}$ is a weight
basis for $\mathfrak{u(\lambda)} \cap \mathfrak{p}$. The vector $u_i \wedge e_{2m+1}$
has weight $\epsilon_i, 1 \leq i \leq r$.
\item We have
$$R(\mathfrak{q}(\lambda)) = i(\lambda) $$
\item Let $W$ be the finite-dimensional  representation with highest weight 
$\lambda$. Then
the representation $A_{\mathfrak{q}(\lambda)}(\lambda)\otimes W^*$ has cohomology in degrees 
$i(\lambda)$ and $n-i(\lambda)$.
\end{enumerate}
\end{lem}
\proof

We will prove only the last statement.
We note that
$L = C \times SO(2s,1)$ with $s = m -i(\lambda)$. Hence $L^c = C \times SO(2s+1)$ and
$L^c/L^c \cap K = S^{2s}$. Thus $A_{\mathfrak{q}}(\lambda)\otimes W^*$ has 
nonzero cohomology in
degrees $i(\lambda)$ and $i(\lambda)+(2(m -i(\lambda))) = 2m -i(\lambda) = n - i(\lambda)$.
\qed

We can now prove Theorem \ref{vanishingtheorem} for the case $n=2m$. Note that
in the statement of the next theorem we are not assuming that
$\mathfrak{q} = \mathfrak{q}(\lambda)$ only that $\mathfrak{q}$ and $\lambda$
are {\em compatible}.

\begin{thm}
$$H^j(A_{\mathfrak{q}}(\lambda) \otimes W_{\lambda}^*) = 0, \ \text{for} \ j< i(\lambda).$$
\end{thm} 
\proof
We claim we have
$$i(\lambda) \leq R(\mathfrak{q}).$$

We deduce the inequality from the equality
$i(\lambda) = R(\mathfrak{q}(\lambda))$ (Lemma \ref{Rlemma})
and the inequality $R(\mathfrak{q}) \geq R(\mathfrak{q}(\lambda))$
from Lemma \ref{reducetostandard}. The claim follows.
Hence $ j < i(\lambda) \Rightarrow j< R(\mathfrak{q}) \Rightarrow
H^j(A_{\mathfrak{q}}(\lambda) \otimes W_{\lambda}^*) = 0$ by Lemma \ref{VZ}.
\qed

\subsection{The case of $SO(2m-1,1)$}

For the case of $SO(2m-1,1)$ we reason as above. Again we have $l=m$.  However 
Lemma \ref{Rlemma} must be modified in case
$i(\lambda) = m$. 

\begin{lem}Suppose that $i(\lambda) = m$. Then there does not
exist $\mathfrak{q}$ such that $\lambda$ is compatible with
$\mathfrak{q}$ and $\lambda|\mathfrak{l}$ is infinitesimally unitary.
\end{lem}

\proof
By assumption $\mathfrak{h} \subset \mathfrak{l}$ for any standard $\mathfrak{q}$.
But we have $\mathfrak{a} \subset \mathfrak{h}$. Now if $\eta = 
(a_1,\cdots,a_m) \in \mathfrak{h}^*$ is infinitesimally unitary then its last component $a_m$
must be pure imaginary. But the last component of $\lambda$ is a {\em nonzero}
integer by assumption. 
\qed

Thus we obtain

\begin{thm}
Let $W$ be a finite-dimensional irreducible representation of
$SO(2m-1,1)$ with highest weight $\lambda$ satisfying $i(\lambda)=m$. 
Then for all $p$ and all cocompact lattices $\Phi \subset SO(n,1)$
we have
$$H^p(\Phi,W) =0.$$
\end{thm}

The case in which $i(\lambda) < m$ is roughly the same as for the case
of $SO(2m,1)$ and is left to the reader.  However we caution the
reader that in this case the space $\mathfrak{u} \cap \mathfrak{p}$
is not $\mathfrak{h}$-stable (because it is not $\mathfrak{a}$-stable),

In particular we have 

\begin{lem}
\begin{enumerate}
\item $L(\lambda) =C \times SO(2s-1,1)$ where $C \subset K$.
\item The set $\{ u_i \wedge e_{2m}: 1\leq i \leq i(\lambda) \}$ is a 
basis (but not a weight basis) for $\mathfrak{u}(\lambda) \cap \mathfrak{p}$.
\item  We have
$$R(\mathfrak{q}(\lambda)) = i(\lambda) $$
\item Let $W$ be the finite-dimensional  representation with highest weight 
$\lambda$. Then
the representation $A_{\mathfrak{q}}(\lambda)\otimes W^*$ has cohomology in degrees 
$i(\lambda)$ and $n-i(\lambda)$.
\end{enumerate}
\end{lem}

Theorem \ref{vanishingtheorem} for the case of $SO(2m-1,1)$ and $i(\lambda) < m$
now follows from the previous lemma in the same way as for the case of $SO(2m,1)$.
  
\section{The intersection theory of cycles with local coefficients }
\subsection{Simplicial homology and cohomology with local coefficients}
\subsubsection{The definition of the groups}
We begin by recalling that the fundamental groupoid of a connected topological 
space $X$
is the category whose objects are the points of $X$ and such that 
for $x,y \in X$ the morphisms,\ $Mor(x,y)$ from $x$ to $y$ are the homotopy classes
of paths from $x$ to $y$.
  
\begin{dfn}
A coefficient system (or local coefficient system or
local system) $\mathcal{G}$ of abelian groups over a topological space $X$
is a covariant functor from the fundamental groupoid of $X$ to the category of abelian 
groups.
\end{dfn}

Thus a coefficient system $\mathcal{G}$ assigns an abelian group 
$\mathcal{G}_x$ to each point $x \in X$ and an isomorphism 
$\tau_{\gamma}$ from $\mathcal{G}_x$ to $\mathcal{G}_y$
to each path $\gamma$ from $x$ to $y$ in such a way so that the same isomorphism
is attached to homotopic (rel end points) paths. 

\begin{ex}
Assume that $X$ is a smooth manifold and $E$ is a flat-bundle over $X$.
Then $E$ induces a local system $\mathcal{E}$ of  vector spaces over $X$.
We will abuse notation henceforth and use $E$ to denote both the flat 
vector bundle and the local system that it induces.
\end{ex}

We now define the homology and cohomology groups of a topological space with
coefficients in a local system  $\mathcal{G}$ of abelian groups over a topogical
space $X$. We will do this assuming that 
 $X$ is the underlying space of a connected simplicial complex $K$. 
 We will define the {\em simplicial} homology and cohomology groups 
with values in $\mathcal{G}$. By the usual subdivision argument one can prove
that the resulting groups are independent of the triangulation $K$. 

We define a $p$-chain with values in 
$\mathcal{G}$ to be a formal sum $\Sigma_{i-1}^m  \sigma_i \otimes c_i$ 
where $\sigma_i$ is an 
{\em ordered} 
$p$-simplex and $c_i$ is an element of the fiber 
of $\mathcal{G}$ over the first vertex of $\sigma_i$. We denote the group of 
such chains by 
$C_p(X,\mathcal{G})$ and define the boundary operator $\partial_p : 
C_p(X,\mathcal{G}) \longrightarrow C_{p-1}(X,\mathcal{G})$
for $\sigma = (v_0,v_1, \cdots, v_p)$ and $c$ an element 
of $\mathcal{G}_{v_0}$ by :
$$\partial_p (\sigma \otimes c) = \sigma_0 \otimes \tau_{(v_0,v_1)}(c) + 
\Sigma_{i=1}^p (-1)^i \sigma_i \otimes c.$$
Here $\tau_{(v_0,v_1)}$ is the isomorphism which is the value of $\mathcal{G}$ on the 
edge $(v_0,v_1)$ and
$\sigma_i,\ 0\leq i \leq p$, is the $i$-th face of $\sigma$. This means that
$\sigma_i = (v_0,v_1,\cdots, \hat{v_i}, \cdots, v_p)$ where $\hat{v_i}$ means that
the vertex $v_i$ is omitted. 
Then $\partial_{p-1} \circ 
\partial_p =0$ and we define the homology groups $H_{*}(X,\mathcal{G})$ of $X$ with coefficients
in $\mathcal{G}$ in the usual way. These groups depend only on the topological space 
$X$  and the local coefficient system $\mathcal{G}$.

In a similar way cohomology groups of $X$ with coefficients in  $\mathcal{G}$ are defined.
A $\mathcal{G}$-valued $p$-cochain on $X$ with values in $\mathcal{G}$ is a function $\alpha$ which assigns to each
ordered $p$-simplex $\sigma $ an element $\alpha(\sigma)$ in  $\mathcal{G}_{v_0}$ where
$v_0$ is the first 
vertex of $\sigma$. The coboundary $\delta_p \alpha$ of a $p$-cochain $\alpha$ is defined on a 
$(p+1)$-cochain $\sigma$ by :
$$\delta_p \alpha(\sigma) = \tau_{(v_0,v_1)}(\alpha(\sigma_0))+
\Sigma_{i=0}^p (-1)^i \alpha(\sigma_i).$$

Then $\delta_{p+1}\circ \delta_p = 0$ and we define the cohomology groups 
$H^{*}(X,\mathcal{G})$
of $X$ with coefficients in $\mathcal{G}$ in the usual way.

In what follows we will need the formulas for $\partial$ and $\delta$
when we express chains and cochains with local coefficients in terms
of flat sections of a flat vector bundle $E$. Note that if $t$ is a 
flat section of $E$ over a face
$\tau$ of a simplex $\sigma$ then it extends to a unique flat section
$e_{\sigma,\tau}(t)$ over $\sigma$. Similarly if we have
a flat section $s$ over $\sigma$ it restricts to a flat section 
$r_{\tau,\sigma}(s)$ over $\tau$. Finally if $\sigma = (v_0,\cdots,v_p)$
we define the $i$-th face $\sigma_i$ by 
$\sigma_i = (v_0,\cdots,\hat{v_i},\cdots,v_p)$. Here $\hat{v_i}$ means
the $i$-th vertex has been omitted.

With these notations we have

$$\partial_p (\sigma \otimes s) = 
\sum _{i=0}^p (-1)^i \sigma_i \otimes r_{\sigma_i,\sigma}(s)$$
and

$$\delta_p (\alpha (\sigma)) = 
\sum _{i=0}^p (-1)^i e_{\sigma,\sigma_i} (\alpha(\sigma_i)).$$

\begin{rem}
If we use flat sections as local coefficients we can use {\em oriented} 
simplices instead of {\em ordered} simplices.
\end{rem}

\subsubsection{Bilinear pairings}
Since our only concern in this paper is with local coefficient systems of vector
spaces we will henceforth restrict to that case (although much of the following
could be carried out in the above generality).
From now on we will regard a simplex with coefficients in $\mathcal{E}$ to be
an oriented simplex $\sigma$ together with a parallel section of
$E|\sigma$ (i.e. a flat lift of $\sigma$) and a cochain with values in
$\mathcal{E}$ to be a rule which assigns to every oriented simplex $\sigma$ a
flat section of the restriction of the bundle $E$ to $\sigma$.

Let $x_0$ be a base--point for $X$. We first define the Kronecker pairing between homology and cohomology with flat vector bundle
coefficients. Let $E,F$and $G$ be flat bundles over $X$.  Assume that 
$\nu: E \otimes F \longrightarrow G$ is a parallel section of $Hom( E\otimes F,G)$.
Let $\alpha$ be a $p$-cochain with coefficients in $E$ and $\sigma \otimes c$ be 
an ordered $p$-simplex with coeffients in $F$ with $\sigma = (v_0,\cdots,v_p)$. 
Then the Kronecker index $<\alpha,c >$ is the element of 
$H_0(X,G)$ defined by:
$$< \alpha,\sigma \otimes c > = \tau_{\gamma}(\nu(\alpha) \otimes c).$$
Here $\gamma$ is any path joining $v_0$ to the base-point $x_0$. The reader will
verify that the right-hand side of the above formula is independent of the choice
of $\gamma$ and that the Kronecker index descends to give a bilinear pairing

$$<\ ,\ >:H^p(X,E) \otimes H_p(X,F) \longrightarrow H_0(X,G).$$

We note that if $G$ is trivial then $H_0(X,G)\cong G_{x_0}$. In particular
we get a pairing 

$$<\ ,\ >:H^p(X,E^*) \otimes H_p(X,E) \longrightarrow \mathbb{R}$$
which is easily seen to be perfect (because the chain groups are vector spaces).
Thus we have an induced isomorphism

$$H^p(X,E^*)\cong H_p(X,E)^* \ \text{or} \ H_p(X,E^*)^* \cong H^p(X,E).$$

The coefficient pairing $\nu:E \otimes F \to G$ also induces cup  products
with local coefficients

$$\cup :H^p(X,E) \otimes H^q(X,F) \longrightarrow H^{p+q}(X,G).$$
and cap products with local coefficients(here we assume $m \geq p$)

$$\cap :H^p(X,E) \otimes H_m (X,F) \longrightarrow H_{m-p}(X,G).$$
These are defined in the usual way using  the ``front-face'' and ``back-face''
of an ordered simplex and pairing the local coefficients using $\nu$.

Finally we will need the following
\begin{thm}
Let $X$ be a compact orientable manifold with fundamental class $[X]$.
Then we have an isomorphism
$$ \mathcal{D}:H^p(X,E)  \longrightarrow H_{n-p}(X,E).$$
given by
$$\mathcal{D}(\sigma) = \sigma \cap [X].$$
\end{thm}

\begin{dfn}
Suppose $[a] \in H_p(X,E)$. We will define the {\em Poincar\'e dual} of
$[a]$ to be denoted $PD([a])$ by
$$PD([a]) = \mathcal{D}^{-1}([a]).$$
\end{dfn}

\subsection{The de Rham theory of cohomology with local coefficients
and the dual of a decomposable cycle}

In this subsection we recall the de Rham representations of the cohomology
groups $H^*(X,E)$ and of the Poincar\'e dual class $PD(Y \otimes s)$.

\begin{dfn}
A differential $p$-form $\omega$ with values in a vector bundle $E$
is an section of the bundle $\bigwedge^p T^*(X) \otimes E$ over $X$.
\end{dfn}

Thus $\omega$ assigns to a $p$-tuple of tangent vectors at $x \in X$ a point
in the fiber of $E$ over $x$. 

Suppose now that $E$ admits a flat connection $\nabla$. 

We can then make
the graded vector space of smooth $E$-differential forms $A^*(X,E)$ into a complex
by defining

\begin{multline}
d_{\nabla}(\omega)(X_1,X_2,\cdots,X_{p+1}) = \sum_{i=1}^p (-1)^{i-1} 
\nabla_{X_i}
(\omega(X_1,\cdots,\widehat{X_i},\cdots,X_{p+1})) \\
+ \sum_{i < j} (-1)^{i+j}
\omega([X_i,X_j],X_1,\cdots,\widehat{X_i},\cdots,\widehat{X_j},\cdots,X_{p+1}) .
\end{multline}

Here $X_i,1\leq i \leq p+1$, is a smooth vector field on $X$.

The following lemma is standard.

\begin{lem}
$$d_{\nabla}^{P+1}\circ d_{\nabla}^p =0.$$
\end{lem}

We now construct a map $\iota$ from $A^p(X,E)$ to the group of
simplicial cochains $C^p(X,E)$ as follows. Let $\omega \in A^p(X,E)$ and
$\sigma$ be a $p$-simplex of $K$. Then  we in a neighborhood $U$ of $\sigma$ we may
write $\omega = \sum_i \omega_i \otimes s_i$ where the $s_i$'s  are parallel sections
of $E|U$ and the $\omega_i$'s are scalar forms. We then define
$$<\iota(\omega),\sigma> = \sum_i (\int_\sigma \omega_i) s_i(v_0).$$

The standard double-complex proof of de Rham's theorem due to Weil, 
see \cite{BottTu}, page 138, yields

\begin{thm}
The integration map $\iota:H^*_{de Rham}(X,E) \longrightarrow H^*(X,E)$
is an isomorphism.
\end{thm}

We will need a special representation in de Rham cohomology with 
coefficients in $E$ of  the cohomology class $PD(Y\otimes s)$. Here we
assume that $Y$ is a compact oriented $p$-dimensional submanifold of $X$
and $s$ is a nonzero parallel section of $E|Y$.

Suppose then that the Poincar\'e dual cohomology class of the homology class
carried by $Y$ is represented in de Rham cohomology of $X$
by the closed $n-p$-form $\omega_Y$ where $\omega_Y$ is supported in a tubular
neighborhood $U$ of $Y$. The parallel section $s$ of $E|Y$ extends to a 
parallel section of $E|U$ again denoted $s$. We extend $\omega_Y \otimes s$ to 
$M$ by making it zero outside of $U$. We continue to use the notation
$\omega_Y \otimes s$ for this extended form. We note that it is standard that
the extension of $\omega_Y$ by zero represents the Poincar\'e dual to
$[Y]$.

\begin{lem}\label{Poincaredual}
Then the de Rham cohomology class Poincar\'e dual to the 
cycle with coefficients $Y \otimes s$
is represented by the bundle--valued form $\omega_Y \otimes s$.
\end{lem}
\proof
We are required to prove that for any $E^*$-valued closed $p$--form $\eta$ we 
have
$$\int_M \eta \wedge \omega_Y \otimes s = \int_{Y \otimes s} \eta = 
\int_Y \<\eta,s\>.$$
But 
$$\int_M \eta \wedge \omega_Y \otimes s = \int_M \<\eta,s\>\wedge \omega.$$

But because $s$ is parallel on $U$ the scalar form $\<\eta,s\>$ is closed
and the lemma follows because $\omega_Y$ is the Poincar\'e dual to $[Y]$.
\qed

\subsubsection{The intersection theory of cycles
with local coefficients}

From now on we assume that $X$ is a compact oriented $n$--manifold and $E,F$ and $G$ are 
flat vector bundles over
$X$ and $\nu: E \otimes F \longrightarrow G$ is a parallel section of 
$Hom( E\otimes F,G)$.
We orient the top-dimensional simplices of $X$ so that their boundaries cancel.
Let $[a]\in H_{p}(X,E)$ and $[b]\in H_{q}(X,F)$. 

\begin{dfn}
The intersection product $[a]\cdot [b] \in H_{p+q -n}(X,G)$ of
$[a]$ and $[b]$ is defined by
$$[a]\cdot [b] = \mathcal{D}(PD([b]) \cup PD([a])),$$
\end{dfn}

Our goal for the rest of this section is to give a geometric
meaning to  the  intersection product. By this we mean we would
like to be able to represent it by a simplicial cycle supported
on the intersection of the simplices comprising $a$ and $b$ where $a$ and $b$ are appropriate representatives
of the classes $[a]$ and $[b]$ and we want this intersection to have
an appropriate structure (i.e. a subcomplex of the correct dimension).
First we need to discuss what it means
for representing simplicial cycles $a$ for $[a]$ and $b$ for $[b]$ to be
in {\em general position}. In what follows if $c$ is a simplicial chain we 
define the carrier $|c|$
to be the subset of $X$ which is the union of the simplices of $c$.

\begin{dfn} 
We will say that a $p$-simplex $a$ is in general position to
a $q$-simplex $b$ if the intersection $a\cap b$ is (a simplex)
of dimension less than or equal to $p+q - n$.

We say the simplicial cycles $a = \sum_i a_i$ ( of dimension $p$) and
$b = \sum_j b_j$ (of dimension $q$)  are in general position if, for all $i,j$,
the simplices $a_i$ and $b_j$ are in general position. Hence the
intersections $a_i \cap b_j$ satisfy 
$$dim(a_i \cap b_j) \leq p+q - n.$$
\end{dfn}

Given two simplicial cycles $a$ and $b$ we can move $|a|$ by an ambient isotopy,
one simplex at a time, so
that each moved simplex, $\overline{a}_i$, is in general position to each 
simplex $b_j$ of $b$ by \cite{Hudson}, Lemma 4.6. Let $|\overline{a}|$ be the 
union of the moved simplices (the moved simplices are not necessarily simplices
in $X$). Thus we have
$h_t: |a| \to X$ with $h_0$ equal to the inclusion of $|a|$ in $X$ and 
the image of $h_1$ equal to $|\overline{a}|$. The proof of \cite{Hudson},
Lemma 4.6 shows that we may choose $h_1$ to be a piecewise linear embedding. 
We may regard
the pair $a$ and $h_1$ as a singular cycle homologous to the original
simplicial cycle $a$ in the singular complex.

Finally, according to \cite{Hudson}, Lemma 1.10,
we may subdivide $a$ and $X$ so that $h_1$ is simplicial.
We use $\widetilde{a}$ to denote the subdivided image of $h_1$.
It is automatic , \cite{Hudson}, Lemma 1.3, that
the resulting subdivision of $b$ is also a subcomplex. 
Hence the intersection $|\widetilde{a}| \cap |b|$ is also a subcomplex.
We will label the simplices in $\widetilde{a}$
as $\widetilde{a}_i$ and label the simplices in the subdivision of $b$ by 
$b_j$ so $i$ and $j$ run through new index sets.
Thus we have arranged
that, for all $i,j$, $\widetilde{a}_i$ and $b_j$ are in general position  
and that $|\widetilde{a}|$,$|b|$ and the intersection $|\widetilde{a}|\cap |b|$ 
are subcomplexes. 
It follows that the 
intersection of $|\widetilde{a}|$ and $|b|$ is a union of simplices
$\{c_{ij} : = \widetilde{a}_i \cap b_j \}$ such that for each pair $i,j$,
$c_{ij}$ is a face common to  $\widetilde{a}_i$ and  $b_j$. We use the simplicial
map $h_1$ to transfer the orientation of each simplex of the subdivision
of $a$ to its image in $\widetilde{a}$ to obtain a simplicial cycle.  
Since $h_1$ is homotopic to the  inclusion of $a$ this cycle is homologous to 
$a$. 

We wish to assign an orientation to each simplex $c_{ij}$ in the
intersection.  If $dim(c_{ij}) < p+q - n$ we assign the coefficient zero to it.
Otherwise we first choose an orientation of $c_{ij}$. 
Since $\widetilde{a}_i$ is oriented and $c_{ij}$ is oriented we obtain an induced
orientation on $Lk(c_{ij},\widetilde{a}_i)$, the link of $c_{ij}$ in $\widetilde{a}_i$
and on $Lk(c_{ij},b_j)$.
Now the star $St(c_{ij})$ is
an $n$-ball and has an orientation induced by the orientation of $X$. 
The union of the vertices of $c_{ij}, 
Lk(c_{ij},\widetilde{a}_i)$ and $Lk(c_{ij},b_j)$ spans 
a subpolyhedron $P$  of $St(c_{ij})$ with nonempty interior . We arrange the 
vertices so that each of the three subsets of vertices induces the correct 
orientation on $c_{ij}$, resp. $Lk(c_{ij},\widetilde{a}_i)$, resp. $Lk(c_{ij},b_j)$. 
We orient $P$ by taking 
the vertices in the order given and assign the coefficient $+1$ to $c_{ij}$with
the chosen orientation if
this orientation agrees with the orientation $P$ receives from the
orientation of $St(c_{ij})$. Otherwise we assign  the coefficient $-1$ 
(or equivalently we take the opposite orientation on $c_{ij}$).

We can now give a geometric definition of  the intersection product of 
simplicial cycles with local coefficients.

\begin{enumerate}
\item Subdivide and choose representing simplicial chains $a = \sum_i a_i \otimes s_i$  and 
$b= \sum_j b_j \otimes t_j$  in general position.

\item Assign to each simplex $c_{ij}$ of the intersection an orientation 
 according to the rule above.

\item Give each $c_{ij}$ as coefficient  the parallel 
section of $G$ over  $c_{ij}$ obtained from the parallel
sections $s_i$ and $t_j$ by restricting each of them to $c_{ij}$
and applying $\nu$.
\end{enumerate}

We will (temporarily) denote the intersection product of $[a]$ and
$[b]$ defined in this way by $[a]\circ[b]$.

\begin{thm}
The geometric definition of the intersection product agrees with the
cup-product definition.
$$[a]\cdot [b] = [a]\circ[b].$$
\end{thm}

We will not prove this theorem but will prove the special case of it we need,
see Theorem \ref{Intersectionformula}.

\subsection{Decomposable cycles and an intersection formula}

There is a particularly simple construction of cycles with coefficients in $E$.
Let $Y$ be a closed, oriented submanifold of $X$ of codimension $p$ and let
$s$ be a parallel section of the restriction of $E$ to $Y$. Let [Y] denote the
fundamental cycle of $Y$ so $[Y]= \Sigma_i \sigma_i$, a sum of ordered simplices.

\begin{dfn}
$Y\otimes s$ denotes the $p$-chain with values in $E$ given by
 $$Y = \Sigma_i \sigma_i \otimes s_i$$
where $s_i$ is the value of $s$ on the first vertex of $\sigma_i$.
\end{dfn}

\begin{lem}
$Y\otimes s$ is an $p$ cycle.
\end{lem}
\proof
Let  $\tau$ be an $p-1$ simplex in $Y$. Then the coefficient of $\tau$ in the 
boundary of $Y\otimes s$ is the product of the coefficient of $\tau$ in the
boundary of the fundamental cycle times the value of $s$ on the first vertex of
$\tau$. But the first factor of the this product is zero.
\qed
We will refer to such cycles henceforth as {\em decomposable} cycles. 
The terminology
is motivated by the case of $1$-cycles and their connection
with the Eilenberg-MacLane complex as we will explain in the next
subsection.

 We now specialize the intersection product of cycles with local coefficients to
cycles of the above type. Suppose that $Y_1$ and $Y_2$ are closed, oriented 
submanifolds of dimensions $p$ and $q$ respectively, $s_1$ and $s_2$ are parallel
sections of $E|Y_1$ and $F|Y_2$ respectively and $Y_1$ and $Y_2$ intersect
transversally. Then we may simplify the previous formula defining the intersection
class $(Y_1 \otimes s_1)\cdot (Y_2 \otimes s_2)$ as follows
\begin{enumerate}
\item Intersect $Y_1$ and $Y_2$ in the usual way to obtain a (possibly
disconnected) oriented dimension $p+q-n$ manifold $Z = \coprod_i Z_i$ with  
intersection multiplicity $\epsilon_i = \pm 1$ along $Z_i$.
\item Assign to $Z_i$ the parallel section of $G|Z_i$ given by 
$\nu(s_1|Z_i,s_2|Z_i)$.
\end{enumerate}
Then 

\begin{thm}\label{Intersectionformula}

$$(Y_1 \otimes s_1)\cdot (Y_2 \otimes s_2) = \sum_i \epsilon_i Z_i \otimes 
\nu(s_1|Z_i,s_2|Z_i)$$

\end{thm}

\proof

By Lemma \ref{Poincaredual} we have
$$ PD(Y_1\otimes s_1) = \omega_{Y_1} \otimes s_1  $$
and 
$$PD(Y_2\otimes s_2) = \omega_{Y_2} \otimes s_2.$$

Now by definition we have
$$PD(Y_1\otimes s_1\cdot Y_2\otimes s_2) = \omega_{Y_2} \otimes s_2 \wedge
\omega_{Y_1} \otimes s_1 = \omega_{Y_2} \wedge 
\omega_{Y_1} \otimes \nu(s_1,s_2).$$ 
But it is well-known, see \cite{Bredon}, Theorem 11.7, that
$PD (\sum_i \epsilon_i Z_i) = \omega_{Y_2} \wedge \omega_{Y_1}$, 
We again apply Lemma \ref{Poincaredual}
to obtain
$$PD(\sum_i \epsilon_i Z_i \otimes 
\nu(s_1|Z_i,s_2|Z_i)) = \omega_{Y_2} \wedge 
\omega_{Y_1} \otimes \nu(s_1,s_2).$$
\qed

\subsubsection{Decomposable cycles and the  Eilenberg-MacLane complex}
In this subsection we will assume that manifolds $X$ and $Y$ are 
Eilenberg-MacLane spaces and will show that decomposable cycles
$Y \otimes s$ have representatives of a particularly simple type in
the Eilenberg-MacLane complex $C_.(\pi_1(X))\otimes V$.
Let $\Gamma$
be the fundamental group of $X$.  The flat bundle $E$ corresponds to a
$\Gamma$-module $V$. The one dimensional simplicial homology group 
$H_1(X,E)$ defined above) is isomorphic to the Eilenberg-Maclane (i.e.group) homology
group $H_1(\Gamma,V)$. Let $\mathbb{R}(\Gamma)$ be the real group ring of $\Gamma$. 
A one-- chain in the Eilenberg-MacLane 
complex of $\Gamma$ with values in $V$ is an element of the tensor product
$\mathbb{R}(\Gamma)\otimes V$; that is, a sum 
$c = \sum_i \gamma_i \otimes v_i$ . The formula for the boundary of the above
Eilenberg-MacLane chain is
$$\partial_1 (c) = \sum_i \gamma_i \cdot v_i - v_i.$$
Thus the simplest kind of Eilenberg-MacLane one-cycle is a decomposable element
$\gamma \otimes v$ where $\gamma\cdot v = v$. Such cycles (decomposable elements of the
group of one chains) were called decomposable cycles in \cite{GoldmanMillson}
(see also \cite{KatokMillson} for the connection with Eichler-Shimura periods
of classical modular forms). Now let 
$c$ be a smooth curve
representing $\gamma$. In case $c$ can be chosen to be embedded then
the vector invariant vector $v$ corresponds to a parallel section $s$ along $c$
and $\gamma \otimes v$ corresponds under the above isomorphism to the cycle $c \otimes s$ 
above. Of course from the point of view of simplical chains such cycles are not 
decomposable elements but they are ``as decomposable as possible" 
whence our terminology.

Assume now that $X$  is a closed oriented $n$-manifold 
which is an Eilenberg-MacLane space of type $K(\Gamma,1)$ 
and $Y$ is a closed oriented $m$-dimensional submanifold 
which is an Eilenberg-MacLane space of type $K(\Phi,1)$
(this will be the case for the cycles studied in this paper).
We will now see that the $m$-cycle $Y\otimes s$ again has a nice interpretation 
in the Eilenberg-MacLane chain complex $C_.(\Gamma) \otimes V$.
Indeed, see \cite{MacLane}, page 114, (5.3), we recall that the Eilenberg-MacLane boundary operator 
$\partial : C_p(\Gamma) \otimes V \to C_{p-1}(\Gamma) \otimes V$ is given
by the formula

\begin{align*}
& \partial((\gamma_1,\gamma_2,\cdots,\gamma_p) \otimes v) = \ 
(\gamma_2,\cdots,\gamma_p) \otimes \gamma_1\cdot v \ \ + \\
& \sum\limits_{i=1}^{p-1}
(-1)^i(\gamma_1,\cdots,\gamma_i\gamma_{i+1},\cdots, \cdots,\gamma_p) \otimes 
v \  + \  (-1)^p (\gamma_1,\cdots,\gamma_{p-1}) \otimes v.
\end{align*}

From this we see that if $[\Phi]$ is the fundamental cycle of $Y$
expressed in the Eilenberg-MacLane complex $C_m(\Phi)$
and $v$ is the $\Phi$-invariant vector in $V$ corresponding to the section $s$
then in the complex $C_m(\Phi)\otimes V$ we have
$$\partial ([\Phi] \otimes v) = (\partial[\Phi]) \otimes v =0.$$ 

Let $f:\Phi \to \Gamma$ be the map on fundamental groups induced
by the inclusion $Y \subset X$. Then $f$ induces chain maps
$\widehat{f} :C_.(\Phi)  \to C_.(\Gamma) $ and
$\widehat{f} \otimes 1 :C_.(\Phi) \otimes V \to C_.(\Gamma) \otimes V$. Indeed,
if $\sigma$ is an Eilenberg-MacLane
$p$-simplex for $\Phi$, we let $\widehat{f}(\sigma)$ denote the simplex (for $\Gamma$) 
obtained by replacing each
group element in $\sigma$ with its image under $f$. We extend $\widehat{f}$ 
linearly to chains.
Then the $m$-cycle with local coefficients $Y\otimes s$ corresponds to the 
Eilenberg-MacLane $m$-cycle $\widehat{f}[\Phi] \otimes v$.

\section{The standard arithmetic quotients of $\mathbb{H}^n$ 
and their decomposable cycles with coefficients}
 
\subsection{The standard arithmetic quotients 
$M= \Gamma \backslash \mathbb{H}^n$}

let $\mathbb{K}$ be a totally-real number field with archimedean completions
$\{v_0,\cdots,v_r\}$. Let $\mathcal{O}$ be the ring of algebraic
integers in $\mathbb{K}$. Let $\underline{V}$ be an oriented vector space over 
$\mathbb{K}$ of dimension $n+1$. Let $V$ be the completion of 
$\underline{V}$ at $v_0$. Let $q:V \to \mathbb{K}$ be a quadratic form 
such that $q$ has signature $(n,1)$ at the completion $v_0$ and is 
positive-definite at all other completions. 
In fact we will restrict ourselves to the example of the Introduction, that is:
$$f(x_1,\cdots,x_{n+1}) = x_1^2 + x_2^2 +\cdots + x_n^2 - \sqrt{m} x_{n+1}^2.$$

Let $\underline{G}$ be
the algebraic group of whose $\mathbb{K}$-points is the group of
orientation preserving isometries of $f$. Let $L$ be an $\mathcal{O}$ 
lattice in $V$ and let $\Phi = \underline{G}(\mathcal{O})$ be the subgroup
of $\underline{G}(\mathbb{K})$ consisting of those elements that take
$L$ into itself. We define $G:=\underline{G}(\mathbb{R})$. We choose
a base-point $o$ in $\mathbb{H}^n$ and let $K$ be the maximal compact
subgroup of $G$ that fixes $o$. We let $\mathfrak{b}$ be an ideal in 
$\mathcal{O}$ and let $\Gamma = \Gamma(\mathfrak{b})$ be the congruence
subgroup of $\Phi$ of level $\mathfrak{b}$ (that is the elements of
$\Phi$ that are congruent to the identity modulo $\mathfrak{b}$).

In what follows it will be convenient to take the projective model for
$\mathbb{H}^n$. We let $D \subset \mathbb{P}(V)$ be the set of 
``negative lines'', that is
$$D= \{z \in \mathbb{P}(V): (\ ,\ )|z \leq 0\}.$$

We define $z_0 = e_{n+1}$ and 
$$V_+ = (\mathbb{R} z_0) ^{\perp}.$$
We take $z_0$ to be the base-point $o$.

\subsection{Totally-geodesic submanifolds in the standard arithmetic
examples}
Let $x=\{x_1,x_2,\cdots,x_k\} \in \underline{V}^k$ (i.e. a $k$-tuple of $\mathbb{K}$
--rational vectors). We assume that $X:= span \ x$ has dimension $k$
and moreover $(\ ,\ )|X$ is positive definite. In what follows we will
let $r_X$ be the isometric involution of $V$ given by

$$r_X(v) =
\begin{cases}
-v &   v\in X \\
v  &  v\in X^{\perp}.
\end{cases}
$$

We define the
totally-geodesic subsymmetric space $D_X$ by
$$D_X = \{\mathbb{R}z \in D : (z,x_i) = 0, 1\leq i \leq p\}.$$
Then $D_X$ is the fixed-point set of $r_X$ acting on $D$.

We also define subgroups $G_X$ (resp. $\Gamma_X$) to be the stabilizer
in $G$ (resp. in $\Gamma$) of the subspace $X$. We define $G^{\prime}_X
\subset G_X$ to be the subgroup that acts trivially on $X$.

We then have, see \cite{JohnsonMillson}, Lemma 7.1 or
\cite{FarrellOtanedaRaghunathan},Lemma 2.4.

\begin{lem}
There exists a congruence subgroup $\Gamma := \Gamma(\mathfrak{b})$ of
$\Phi$ such that

\begin{enumerate}
\item $M :=\Gamma\backslash D$ is an orientable
compact manifold of codimension $k$.
\item
The image $M_X$ of $D_X$ in $M$ is the quotient $\Gamma_X\backslash D_X$
and is an orientable submanifold.
\end{enumerate}
\end{lem}

In what follows we will require $x \in V_+^k$. Accordingly we
have $z_0 \in D_X$.

\subsection{Attaching a local coefficient $s_x$ to the cycle $M_X$}
\subsubsection{The existence of a parallel section}
We now want to promote $M_X$ to a (decomposable) cycle with coefficients
for appropriate coefficient systems $W$ by finding a nonzero 
parallel section of $s_X$ of $\widetilde{W}|M_X$.
The main point is

\begin{lem}
If $\Gamma$ is a neat subgroup (\cite{Borel},pg.117) then $\Gamma_X$
acts trivially on $X$.
\end{lem}
\proof
We have a projection map $p_X :\Gamma_X \to O(X_1) \times O(X_2) \times \cdots
\times O(X_r)$. Here by $X_i$ we mean the $i$-th completion of $X$.
The $i$-th completion of $(\ ,\ )$ restricted to $X_i$ is positive
definite for $1\leq i \leq r$. Furthermore the splitting $V=X \oplus
X^{\perp}$ is defined over $\mathbb{K}$. Thus the diagonal embedding
of the intersection $L_X = L \cap X$ is a lattice in $\oplus_{i=0}^r 
X_i$ which is invariant under $p_X(\Gamma_X)$. Hence $p_X(\Gamma_X)$
is a discrete subgroup of a compact group hence a finite group. Hence
if $\gamma \in p_X(\Gamma_X)$ then all eigenvalues of $\gamma$ are roots
of unity. Since $\Gamma$ is neat all eigenvalues must be $1$ and the
lemma follows.
\qed

We can now prove the existence of a nonzero parallel section along $M_X$
provided that $dim(X)$ is big enough. For economy of notation we will
adopt the following terminology up to the end of the following
proposition. We will say a dominant weight $\mu$ for $SO(2m-1,1)$
is {\em admissible} if $i(\mu) \leq m-1$. For the case $SO(2m,1)$
all dominant weights will be defined to be admissible. We note by Theorem \ref{vanishingtheorem} 
that if the highest weight of $W$ is not admissible then all cohomology
groups with coefficients in $W$ vanish.

\begin{prop}\label{existenceofinvariant}
Suppose that $\Gamma$ is neat and $\mu$ is admissible.  Then there exists a 
nonzero vector $\tau \in S_{[\mu]}V \ \text{invariant under} \  \Gamma_X 
\Leftrightarrow dim(X) \geq i(\mu)$.
\end{prop}
\proof
We apply the Gelfand-Tsetlin branching theorem for the pair of groups
$SO(m-1) \subset SO(m)$, \cite{Boerner}, pg. 267 and pg. 269 to
find that if the irreducible representation of $SO(m-1)$ with highest weight
$\nu$ occurs in the restriction to $SO(m-1)$ of the representation of
$SO(m)$ with admissible highest weight $\mu$ then 
$$i(\mu) -1 \leq i(\nu) \leq i(\mu) \ \text{and}\  i(\nu) = i(\mu) -1 
\ \text{is realized}.$$

Thus we have to branch at least $i(\mu)$ times to get the zero weight
and if we branch $i(\mu)$ times then the zero weight occurs.
\qed

In what follows we will need an explicit formula for an invariant vector
because we will have to take the inner product of two such invariants.
The rest of this section is devoted to finding such an explicit formula.
However the result that in order to find a parallel section it is
necessary that
$$dim(X) \geq i(\mu)$$
will play a critical role throughout the rest of the paper.

We will use Weyl's construction of the irreducible representations
of $O(V)$ to find an explicit formula for a parallel section.
We first review Weyl's construction.

\subsubsection{The harmonic Schur functors}

We will follow \cite{FultonHarris} in our description of the harmonic
Schur functor $U \to S_{[\mu]}U$ corresponding to a partition $\mu$
on quadratic spaces $U, (\ ,\ )$. Here $(\ ,\ )$ is assumed to
be nondegenerate. There is also a nice treatment in 
\cite{GoodmanWallach}. We assume that $U$ has dimension $m$.

Suppose then that $W_{\mu}$ is the irreducible representation of $SO(U)$ 
with
highest weight $\mu = (b_1,b_2,\cdots,b_l)$ where $l= [\frac{m}{2}]$.
We will abuse notation and use $\mu$ to denote the corresponding partition
of $d = \sum b_i$. We extend the quadratic form $(\ ,\ )$  to $\otimes ^d U$
as the $d$-fold tensor product and note that the action of $S_d$ on
$\otimes ^d U$ is by isometries. For each pair $I = (i,j)$ of integers
between $1$ and $d$ we have the contraction operator $\Phi_I
:\otimes^d U \to \otimes^{d-2} U$, the operator 
$\Psi_I:\otimes^{d-2} U \to \otimes^d U$ that inserts the  (dual of
the form ) form $(\ ,\ )$ an into the $(i,j)$-th spots and the composition
$\theta_I := \Psi_i \circ \Phi_I$, (the notation is that of 
\cite{FultonHarris}). We define the harmonic $d$-tensors, to be
denoted $U^{[d]}$, to be the kernel of all the contractions $\Phi_I$.
Following \cite{FultonHarris}, pg. 263, we define the subspace 
$U^{[d]}_{d-2r}$ of $U^{\otimes d}$ by

$$U^{[d]}_{d-2r} = \sum \Psi_{I_1} \circ \cdots \Psi_{I_r}U^{[d -2r]}.$$

Carrying over the proof of \cite{FultonHarris}, Lemma 17.15 (and the
exercise that follows it) from the symplectic case to the orthogonal case
we have

\begin{lem}
We have a direct sum, orthogonal for $(\ ,\ )$,
$$U^{\otimes d} = U^{[d]} \oplus \oplus_{r=1}^{[\frac{d}{2}]} 
U^{[d]}_{d - 2r}.$$
\end{lem}

We define the harmonic projection $\mathcal{H}:\otimes V \to \otimes V$ to be the 
orthogonal projection on the harmonic $d$-tensors $U^{[d]}$. The space of
harmonic $d$-tensors $U^{[d]}$ is invariant under the action of
$S_d$. Consequently we may apply the idempotents in the group algebra
of $S_d$ corresponding to partitions to further decompose $U^{[d]}$
as an $SO(U)$--module.

We use $\mu$ to denote the associated partition of $d$. {\em We note that
$i(\mu)$ equals the length of the first column of $\mu$.} We label the Young 
diagram corresponding to the partition $\mu$ in the standard way, i.e.
by putting $1,2,\cdots,b_1$ in the first row
and continuing in the obvious way to obtain a standard tableau $T$.
Let $P$ (resp. $Q$) be the group preserving the rows (resp. 
columns) of $T$. Define elements of the group ring of $S_d$
by  $\mathcal{P} = \sum_P p$ and $\mathcal{Q} = \sum_Q \epsilon(q) q$.
Again following \cite{FultonHarris}, pg. 296, we define the harmonic Schur 
functor $S_{[\mu]}U$ as follows.

\begin{dfn}
$$S_{[\mu]}U= \mathcal{Q}\mathcal{P} U^ {[d]}.$$
\end{dfn}

We then have the following theorem, \cite{FultonHarris}, Theorem 19.19

\begin{thm}
The space $S_{[\mu]}U$ is an irreducible representation of $O(V)$.
It is nonzero if and only if the sum of the first two columns of
the partition $\mu$ is less than or equal to $dim(U)$.
\end{thm}

In what follows we will need the following

\begin{lem}
\begin{enumerate}
\item $\mathcal{H},\mathcal{P}$ and $\mathcal{Q}$ are self-adjoint
relative to $(\ ,\ )$.
\item $\mathcal{H}$ commutes with $\mathcal{P}$ and $\mathcal{Q}$.
\end{enumerate}
\end{lem}

\proof
It is clear that $\mathcal{H}$ is self-adjoint. The arguments for $\mathcal{P}$ and
$\mathcal{Q}$ are the same. We give the one for $\mathcal{Q}$.
We will use the symbol $q$ to denote both the element $q \in Q$
and the corresponding operator on $V^{\otimes d}$. Since $q$
is an isometry we have $q^* = q^{-1}$. Hence we have
$$\mathcal{Q}^* = \sum \epsilon(q)q^* = \sum \epsilon(q^{-1}) q^{-1}
= \mathcal{Q}.$$

To prove that $\mathcal{H}$ commutes with $\mathcal{P}$ and $\mathcal{Q}$
it suffices to prove that $\mathcal{H}$ commutes with every element 
$g \in S_d$. But $S_d$ acts by isometries and preserves $U^{[d]}$.
Consequently it commutes with orthogonal projection on $U^{[d]}$.
\qed

\subsubsection{An explicit formula for a parallel section}
We will now return to the case in hand.
Assume  that $i(\mu) \leq k \leq m = [\frac{n+1}{2}]$.
We now construct an explicit parallel section $s_e$
where $e=(e_1,e_2,\cdots,e_k)$. Suppose that $b_i =0, i\geq k+1$.
Put $E_k = span \ \{e_1,\cdots,e_k\}.$
Define an element $\sigma_e \in S^{b_1}E_k \otimes S^{b_2}E_k
\otimes \cdots S^{b_k}E_k$ by
$$ \sigma_e = e_1^{\otimes b_1}\otimes \cdots e_k^{\otimes b_k}.$$
Clearly we have $\mathcal{P}\sigma_e = \sigma_e$.

We define $\tau_e \in S_{[\mu]}V$ by 
$$\tau_e = \mathcal{Q}\mathcal{P}\mathcal{H} \sigma_e.$$
We observe that $\mathcal{H}\sigma_e \notin \otimes E_k$
and $\tau_e \notin S_{[\mu]}E_k$ because the harmonic projection $\mathcal{H}$ 
corresponding to the pair $V,(\ ,\ )$ does not carry
$S^*(E_k)$ into itself.
We have

\begin{lem}
$\tau_e$ is $G_{E_k}^{\prime}$--invariant.
\end{lem}
\proof
$g \cdot \mathcal{Q}\mathcal{P}\mathcal{H} \sigma_e =
\mathcal{Q}\mathcal{P}\mathcal{H} g \cdot \sigma_e 
=\mathcal{Q}\mathcal{P}\mathcal{H} \sigma_e.$
\qed

Continue to assume that $i(\mu) \leq k \leq m =[\frac{n+1}{2}]$.
Recall that we previously defined elements $u_i, v_i \in V, 1\leq i \leq 
m=[\frac{n+1}{2}]$ by 
$u_i = e_i - \imath e_{m+i}, v_i = e_i+ \imath e_{m+i}, 1 \leq i \leq m-1$.
Furthermore we have
 $$u_m = \begin{cases} 
e_m - \imath e_{2m}, n = 2m\\
e_m - e_{2m}, n = 2m -1
\end{cases}$$
and 
$$v_m = \begin{cases} 
e_m + \imath e_{2m}, n = 2m\\
e_m + e_{2m}, n = 2m -1.
\end{cases}$$ 

We next define the element  $\sigma_u \in S_{[\mu]}V \otimes \mathbb{C}$
by
$$ \sigma_u = u_1^{\otimes b_1}\otimes \cdots u_k^{\otimes b_k}.$$

We recall that $(u_i,u_j)=0,1\leq i,j \leq k$ and we obtain

\begin{lem}
$\sigma_u$ is a harmonic tensor.
\end{lem}

\begin{cor}
$\tau_u = \mathcal{Q}\sigma_u = \mathcal{Q} \mathcal{P} \mathcal{H}\sigma_u 
\in S_{[\mu]}X \otimes \mathbb{C}$.
\end{cor}

\begin{rem}
The vector $\tau_u$ is a weight vector of weight $(b_1,\cdots,b_m)$
and consequently is a highest weight vector for $W$.
\end{rem}
We can now prove that $\tau_e$ and $\tau_u$ are both nonzero.

\begin{lem}\label{nonvanishinginvariant}
$(\tau_e,\tau_u) = 1$.
\end{lem}

\proof
$$(\tau_e,\tau_u) = (\mathcal{Q}\mathcal{P}\mathcal{H}\sigma_e,\mathcal{Q}\sigma_u)
=(\mathcal{P}\mathcal{H}\sigma_e,\mathcal{Q}\sigma_u)= 
(\mathcal{H}\mathcal{P}\sigma_e,\mathcal{Q}\sigma_u) =
(\sigma_e,\mathcal{H}\mathcal{Q}\sigma_u) = (\sigma_e,\mathcal{Q}\sigma_u).$$
Next we note that 
$(e_i,u_j) = \delta_{ij}$.
Hence for $q \in Q_{\mu}$ we have
$$(\sigma_e,q\sigma_u)=\begin{cases} & 1 , \  q=1 \\
& 0, \  q \neq 1.
\end{cases}$$
\qed

We can now give a formula for a nonzero invariant. We have proved

\begin{thm}\label{explicitformula}
Suppose that  $X\subset V_+$ satisfies 
$i(\mu) \leq dim(X) = k \leq m=[\frac{n+1}{2}]$.
Assume that $\mu = (b_1,\cdots, b_k,0,\cdots,0)$.
Suppose that  $x = (x_1,\cdots,x_k)$ is an orthogonal basis for $X$. Define
$$\tau_x = \mathcal{Q} \mathcal{P} \mathcal{H} ( x_1^{b_i} \otimes
\cdots \otimes x_i^{b_k}).$$
Then $\tau_x$ is
an nonzero $G_X^{\prime}$--invariant in $S_{[\mu]}V$.
\end{thm}
 
 Finally we will need to find $k$-tuples $y = (y_1,\cdots,y_k) \in V^k$
 such that if we define $ \sigma_y = y_1^{\otimes b_1}\otimes 
 \cdots y_k^{\otimes b_k}$ and 
 $\tau_y = \mathcal{Q}\mathcal{P}\mathcal{H}\sigma_y$ we have
 
 \begin{lem}\label{innerproduct}
 $$(\tau_x,\tau_y) \neq 0.$$
 \end{lem}

\proof
There exist $s_1 \in S^{b_1}V,\cdots, s_k \in S^{b_k}V$ such that
$(\tau_x,\mathcal{Q}\mathcal{P}\mathcal{H}(s_1 \otimes \cdots s_k)) \neq 0$.
But the pure powers $v^m, v \in V$ span $S^mV$. Hence, we conclude that
there exist $y_1,\cdots y_k \in V$ such that
$$(\tau_x, \mathcal{Q}\mathcal{P}\mathcal{H}(y_1^{\otimes b_1}\otimes \cdots 
\otimes
y_k^{\otimes b_k})) \neq 0.$$
\qed

\section{Configurations of compact totally-geodesic submanifolds
in the standard arithmetic examples}

In this section we will review the results of \cite{JohnsonMillson}
and \cite{FarrellOtanedaRaghunathan}. 

Let $x= \{x_1,x_2,\cdots,x_k\}$ and $y=\{y_1,y_2,\cdots,y_l\}$ be
such that $z = x \cup y$ spans a subspace $Z$ of $V$ of dimension $k+l$
such that the restriction of $(\ ,\ )$ to $Z$ is positive definite.
In what follows we let
$X = span\  x$ and $Y= span\  y$.

We have the Theorem 7.2 of \cite{JohnsonMillson}, see also Corollary 2.26
of \cite{FarrellOtanedaRaghunathan}.

\begin{thm}\label{onecomponent}
There exists a congruence cover $M^{\prime}$ of $M$ such that 
$M_X^{\prime}\cap M_Y^{\prime}$ consists of the single component
$M_Z^{\prime}:= \pi^{\prime}(D_Z)$. Moreover for
any congruence cover $M^{\prime \prime }$ of $M^{\prime}$ the
intersection $M_X^{\prime \prime}\cap M_Y^{\prime \prime}$
again consists of the single component 
$M_Z^{\prime \prime}:= \pi^{\prime \prime}(D_Z)$
\end{thm}

\begin{cor}
Suppose $k+l=n$. Then there exists a congruence cover $M^{\prime}$ of $M$ such that 
$M_X^{\prime}\cap M_Y^{\prime}$ consists of the single point $\pi^{\prime}(z_0)$.
\end{cor}

\section{The nonvanishing theorem}

We can now prove the desired nonvanishing theorem.
Suppose $W$ is an irreducible representation of $SO(n,1)$ with highest
weight $\mu = (b_1,\cdots,b_m)$. Here we assume that if $n$ is even then
$n=2m$ and if $n$ is odd then $n=2m-1$. Assume $i(\mu) \leq k \leq n$ and choose
$x=e=(e_1,\cdots,e_k)$. We have shown that there exists a nonzero parallel
section $s_e$ of $\widetilde{W}|M_X$ corresponding to the
$\Gamma_X$-invariant vector $\tau_e \in W$. But we will also need a
parallel section over a complementary cycle $M_Y$. Hence by 
Proposition \ref{existenceofinvariant} we require
$$ k \geq i(\mu) \ \text{and} \ n-k \geq i(\mu).$$ 
Hence we require $i(\mu) \leq [\frac{n}{2}]$. We will also
need to construct the section $s_y$ from a $k$-element subset of a spanning 
set $y$ for $Y$.
Accordingly we will assume henceforth that $k \leq n-k$ or 
$k \leq [\frac{n}{2}]$.

\begin{lem}\label{tuples}
Assume now that $i(\mu) \leq k \leq [\frac{n}{2}]$ whence $k \leq n-k$.
Then there exist (infinitely many) $\mathbb{K}$-rational $n-k$-tuples 
$y=(y_1, \cdots,y_{n-k}) 
\in V^{n-k}$ such that if we define $y^{\prime} = (y_1, \cdots,y_k)$ then we have
\begin{enumerate}

\item $(\tau_e,\tau_{y^{\prime}}) \neq 0$
\item Let $z = \{x_1, \cdots, x_k, y_1,\cdots,y_{n-k}\}$ and $Z = span \  z$.
Then $dim(Z) = n$.
\item The restriction of $(\ ,\ )$ to $Z$ is positive definite.
\end{enumerate}

\end{lem}
\proof
We claim that the  negation of each of the first two conditions  defines a 
{\em proper} algebraic subvariety of  $V^{n-k}$. For the second  
the claim is
obvious. The first requires some more work.
 
Define 
$\Phi:V^{n-k} \to S_{[\mu]}V$ by
$$\Phi(y) = \tau_{y^{\prime}}: = \mathcal{Q} \mathcal{P} \mathcal{H} 
(y_1^{\otimes b_i} 
\otimes \cdots \otimes y_k^{\otimes b_k}).$$
Then $\Phi$ is a polynomial mapping.
We define a polynomial $P_e$ on $V^{n-k}$ by
$$P_e(y) = (\tau_e,\tau_{y^{\prime}}).$$ 
But in Lemma \ref{innerproduct} we have proved that there exists
$y \in V^{n-k}$ such that 
$P_e((y_1,\cdots,y_{n-k})) = (\tau_e,\tau_{y^{\prime}}) \neq 0$
Thus $P_e$ is not identically zero and the zero set
of $P_e$ (the negation of 1.) is a proper algebraic hypersurface.

Thus the set of $y$'s  satisfying each of the first two conditions is
nonempty,open and dense.
Finally the set of $y's$ satisfying the third condition is open (in the
classical topology) and consequently the set of points satisfying
all three conditions is a nonempty open subset of $V^{n-k}$. Hence, it 
contains infinitely many rational points.
 
\qed

We can now prove a preliminary nonvanishing result. 

\begin{lem}
Assume $i(\mu) \leq k \leq [\frac{n}{2}]$.
Then $H^{k}(\Gamma,W) \neq 0$.
\end{lem}

\proof
Let $X = E_k$ with $k \leq [\frac{n}{2}]$.
Choose $\{y_1,\cdots,y_{n-k}\} \subset V^{n-k}$ satisfying the three
conditions of Lemma \ref{tuples}
of the previous lemma. Now apply Theorem \ref{onecomponent} and
choose a neat congruence subgroup of $\Gamma$  so that 
$M_X \cap M_Y$ consists of single point. Then 
 
$$M_X\otimes s_x \cdot M_Y \otimes s_y = (s_x,s_y) = (\tau_e,\tau_y) \neq 0.$$
Since $W$ is self-dual we have $H_k(\Gamma,W)^* = H^k(\Gamma,W) \neq 0$.
\qed

\begin{rem}
We see that the method fails (as it should) for the case of
weights $\mu$ with $i(\mu) = m$ for $SO(2m-1,1)$. Indeed in this case we
need both the cardinalities of $x$ and $y$ to be at least $m$ (in order  that the
coefficients $\tau_x$ and $\tau_y$  in $S_{[\mu]}V$ exist). 
But $(\ ,\ )|(X+Y)$ positive definite implies $dim(X+Y) \leq 2m-1$ and 
consequently $X \cap Y \neq 0$.
\end{rem}

Now we have the required nonvanishing theorem.
\begin{thm}\label{nonvanishing}
\begin{enumerate}
\item Suppose that $n$ is even. Then
$$H^k(\Gamma,W) \neq 0, \ i(\mu) \leq k \leq n - i(\mu).$$
\item Suppose $n$ is odd and $i(\mu) < \frac{n+1}{2}$. Then
$$H^k(\Gamma,W) \neq 0, \  i(\mu) \leq k \leq n - i(\mu).$$
\end{enumerate}
\end{thm}
\proof
By Poincar\'e duality
$$H^k(\Gamma,W) \cong H^{n-k}(\Gamma,W).$$
But by the previous lemma we have shown that we have nonvanishing up to
and including the middle dimension in case $n$ is even and up to and
including the first of the two middle dimensions in case $n$ is odd.
\qed

\section{Nonvanishing of cup-products}

We have
\begin{thm}
Let $W_1$ and $W_2$ be irreducible representations with highest weights
$\mu_1$ and $\mu_2$ satisfying $i(\mu_1) = p_1$ and $i(\mu_2) = p_2$.
Suppose that $q_1 \geq p_1$ and $q_2 \geq p_2$ and $q_1 + q_2 \leq 
[\frac{n}{2}]$.
Then (for $\mathfrak{b}$ depending on $W_1$ and $W_2$) the cup-product
$$H^{q_1}(\Gamma,W_1) \otimes H^{q_2}(\Gamma,W_1) \to 
H^{q_1 + q_2}(\Gamma,W_1 \otimes W_2)$$
is a nonzero map. Furthermore if $\pi:W_1 \otimes W_2 \to W$ is a
homomorphism of $SO(n,1)$-modules such that $\pi(\tau_{x^{\prime}} \otimes
\tau_{x^{\prime}}) \neq 0$ then the induced cup-product with values in
$W$ is also nonzero.
\end{thm}

The proof will follow the same lines as the proof of Theorem \ref{nonvanishing}.

Choose $x = (e_1,\cdots,e_{q_1})$ and 
$y= (e_{q_1 + 1},\cdots,e_{q_1+q_2})$. We put 
$x^{\prime}= (e_1,\cdots,e_{p_1})$ and $y^{\prime} = 
(e_{q_1 +1},\cdots,e_{q_1 + p_2})$ and
$z=(e_1,\cdots,e_{q_1}, e_{q_1 +1},\cdots,e_{q_1 + q_2})$.
We define $\tau_{x^{\prime}} \in W_1$ and $\tau_{y^{\prime}} \in W_2$ as before.
We obtain corresponding parallel sections $s_{x^{\prime}}$ and
$s_{y^{\prime}}$ of $\widetilde{W_1}|M_X$ (resp. $\widetilde{W_2}|M_Y$).

We also have the totally geodesic submanifold $M_Z$ where
$Z = span \ z$. After passing to a sufficiently deep congruence subgroup
we have by Theorem \ref{onecomponent}
$$M_X \otimes s_{x^{\prime}} \cdot M_Y \otimes s_{y^{\prime}}
= M_Z \otimes (s_{x^{\prime}} \otimes s_{y^{\prime}}).$$

Here we have used the natural map $\phi$ from $\Gamma(M_Z,\widetilde{W_1}|M_Z) \otimes
\Gamma(M_Z,\widetilde{W_2}|M_Z)$ into 
$\Gamma(M_Z,\widetilde{W_1}|M_Z \otimes \widetilde{W_2}|M_Z)$ given by
$\phi(s_1 \otimes s_2)(x) = s_1(x) \otimes s_2(x)$.
We note if $s_1$ and $s_2$ are parallel then so is $\phi(s_1 \otimes s_2)$.

Thus to prove the theorem we have to prove that the cycle on
the right-hand side is not a boundary. For this we need an
analogue of Lemma \ref{tuples}. The proof is the same as before.

\begin{lem}
Assume that $p_1 \leq q_1$, $p_2 \leq q_2$ and $q_1 + q_2 \leq n - (q_1 + q_2)$.
Then there exist (infinitely many) $\mathbb{K}$-rational $n-k$-tuples 
$$w=(w_1, \cdots,w_{q_1},w_{q_1 +1}, \cdots,w_{q_1 +q_2}, \cdots ,
w_{n-(q_1 +q_2)})$$ 
in $V^{n-(q_1 +q_2)}$ such that if we define $w^{\prime} = 
(w_1, \cdots,w_{p_1})$ and $w^{\prime \prime} = 
(w_{q_1 + 1}, \cdots,w_{q_1 + 1 + p_2})$ then we have
\begin{enumerate}

\item $(\tau_{x^{\prime}},\tau_{w^{\prime}}) \neq 0$.
\item $(\tau_{x^{\prime \prime}}, \tau_{w^{\prime \prime}}) \neq 0$.
\item Let $u = \{e_1,\cdots, e_{q_1 + q_2},w_1,  
\cdots, w_{n-(q_1 +q_2)}\}$ and $U = span \  u$. Then $dim(U) = n$.
\item The restriction of $(\ ,\ )$ to $U$ is positive definite.
\end{enumerate}

\end{lem}

We can now prove the Theorem. Put $W = span \  w$. Form the decomposable cycle
$M_W \otimes (s_{w^{\prime}} \otimes s_{w^{\prime \prime}})$ and 
intersect with $M_Z \otimes (s_{x^{\prime}} \otimes s_{y^{\prime}})$
(after passing to a deep congruence subgroup so that
$M_Z \cap M_W$ is a point).
We obtain
$$M_Z \otimes (s_{x^{\prime}} \otimes s_{y^{\prime}}) \cdot 
M_W \otimes (s_{w^{\prime}} \otimes s_{w^{\prime \prime}}) = 
(s_{x^{\prime}}, s_{w^{\prime}}) (s_{y^{\prime}}, s_{w^{\prime \prime}})
=(\tau_{x^{\prime}},\tau_{w^{\prime}}) 
( \tau_{y^{\prime}}, \tau_{w^{\prime \prime}}) \neq 0.$$

We conclude with an example that shows that $p_1 + p_2 \leq [\frac{n}{2}]$
is not a strong enough assumption to guarantee nonvanishing of the
above cup-product.

\begin{ex}
We take $G=SO(6,1), W_1=\R$ and $W_2 = \bigwedge^3V$ where $V$ is the
standard representation of $G$. We consider the cup-product

$$H^{1}(\Gamma,\R) \otimes H^{3}(\Gamma,\bigwedge^3V) \to 
H^{4}(\Gamma,\bigwedge^3V).$$

But by Poincar\'e duality $H^{4}(\Gamma,\bigwedge^3V) \cong 
H^{2}(\Gamma,\bigwedge^3V) =0$. Note that $p_1 + p_2 = 3 = [\frac{n}{2}]$
whereas $q_1 + q_2 = 4 > [\frac{n}{2}]$.

\end{ex}


\begin{thebibliography}{BaBE}

\bibitem[Bo]{Boerner}
H.\ Boerner, 
``Representations of Groups with Special Consideration for the 
Needs of Modern Physics'', North Holland, 1970.


\bibitem[B]{Borel}
A.\ Borel, 
``Introduction aux Groupes Arithmetiques'', Hermann, 1969.

\bibitem[BLS]{BorelLabesseSchwermer}
A.\ Borel, J.\ P.\ Labesse and J.\ Schwermer,
{\em On the cuspidal cohomology of $S$-arithmetic groups over number fields},  
Compos. Math., vol. {\bf  102 } (1966), pg. 1--40.

\bibitem[BS]{BorelSerre}
A.\ Borel and J.\ P.\ Serre,
{\em Corners and arithmetic groups},  
Comment. Math. Helv., vol. {\bf  48 } (1973), pg. 436--491.


\bibitem[BW]{BorelWallach}
A.\ Borel and N.\ Wallach, 
``Continuous Cohomology, Discrete Subgroups and Representations of
Reductive Groups'', Annals of Math. Studies, vol.  {\bf 94}, Princeton
University Press, 1980.

\bibitem[BT]{BottTu}
R.\ Bott and L.\ W.\ Tu,
`` Differential Forms in Algebraic Topology'', Graduate Texts in Mathematics,
vol.  {\bf 82}, Springer,1982.

\bibitem[Br]{Bredon}
G.\ E.\ Bredon, `` Topology and Geometry'', Graduate Texts in Mathematics,
vol.  {\bf 139}, Springer,1993.

\bibitem[C]{Clozel}
L.\ Clozel,
{\em Sur une question d'Armand Borel}, C.R. Acad. Sci. Paris, S\'er. I 
 vol, {\bf 324}, no. 9, pg. 973--976.



\bibitem[FOR]{FarrellOtanedaRaghunathan}
F.\ T.\  Farrell, P.\ Otaneda and M.\ S.\  Raghunathan,
{\em Nonunivalent harmonic maps homotopic to diffeomorphisms},  
Journ. Diff. Geom. , vol. {\bf 54} (2000), no. 2, pg. 227--253.

\bibitem[FH]{FultonHarris}
W.\ Fulton and J.\ Harris,
``Representation Theory, A First Course'', Graduate Texts in Mathematics,
no. {\bf 129}, Springer, 1991.

\bibitem[FM1]{FunkeMillson1}
J.\ Funke and J.\ J.\ Millson, {\em Cycles in hyperbolic manifolds of
non-compact type and Fourier coefficients of Siegel modular forms},
Manuscripta Math, vol. {\bf 107}, pg. 409--449.


\bibitem[FM2]{FunkeMillson2}
J.\ Funke and J.\ J.\ Millson, in preparation.

\bibitem[GM]{GoldmanMillson}
W.\ M.\ Goldman and J.\ J.\ Millson,
{\em Eichler-Shimura homology and the finite generation of cusp forms
by hyperbolic Poincar\'e series}, Duke Math. J., vol. {\bf 53} (1986),
pg. 1081--1091.

\bibitem[GW]{GoodmanWallach}
R.\ Goodman and N.\ R.\ Wallach,
``Representations and Invariants of the Classical Groups'', 
Encyclopedia of Mathematics and its Applications, 
no. {\bf 68}, Cambridge University Press, 1998.


\bibitem[Ha1]{Harder1973}
G.\ Harder,
{\em On the cohomology of $SL_2(\mathcal{O})$}, Proc. of the Summer School
of Group Representations, I. M. Gelfand ed., Hilger, 1975, London, pg. 139
--150.

\bibitem[Ha2]{Harder1975}
G.\ Harder,
{\em On the cohomology of discrete arithmetically defined groups}, 
Proc. Int. Colloq. on Discrete Subgroups and Applications to Moduli
Problems, Bombay, 1973, Oxford University Press, 1975, pg 129--160.


\bibitem[Hu]{Hudson}
J.\ F.\ P.\ Hudson,
``Piecewise Linear Topology'', Mathematics Lecture Note Series,
Benjamin, 1968.


\bibitem[JM]{JohnsonMillson}
D.\ Johnson and J.\ J.\ Millson,
{\em Deformation spaces associated to compact hyperbolic manifolds},
in ``Discrete Groups in Geometry and Analysis, Papers in Honor of
G.D.Mostow on His Sixtieth Birthday'', Progress in Mathematics vol. {\bf 67},
(1987),pg. 48--106.

\bibitem[KaM]{KatokMillson}
S.\ Katok and J.\ J.\ Millson,
{\em Eichler Shimura homology, intersection numbers and rational structures
on spaces of modular forms},
T.A.M.S., vol. {\bf 300} (1987), pg. 737-757.


\bibitem[KuM]{KudlaMillson}
S.\ S.\ Kudla and J.\ J.\ Millson,
{\em Intersection numbers of cycles on locally symmetric spaces and Fourier
coefficients of holomorphic modular forms in several complex variables},
Publ. Math. I.H.E.S, vol. {\bf 71} (1990), pg. 121-172.


\bibitem[Li]{Li}
J.-S.\ Li,
{\em Nonvanishing theorems for the cohomology of certain arithmetic
quotients},
J. f\"ur die reine und angewandte Mathemtik, vol. {\bf 428} (1992), pg. 177-217.

\bibitem[LS]{LiSchwermer}
J.-S.\ Li and J.\ Schwermer,
{\em Automorphic representations and cohomology of arithmetic groups},
To appear in : Challenges for the $21$-st Century , Proc. Int. Conf. on 
Fundamental Sciences, Singapore, 2001, pg. 102--137.

\bibitem[McL]{MacLane}
S.\ MacLane,
``Homology'', Grundlehren der Mathematischen Wissenschaften, Band {\bf 114},
Springer, 1963.



\bibitem[M1]{Millson1976}
J.\ J.\ Millson,
{\em On the first Betti number of a constant negatively curved manifold},
Ann. of Math.,vol. {\bf 104} (1976), pg. 235--247.

\bibitem[M2]{Millson1985}
J.\ J.\ Millson,
{\em A remark on Raghunathan's vanishing theorem},
Topology,vol. {\bf 24} (1985), pg. 495--498.




\bibitem[MR]{MillsonRaghunathan}
J.\ J.\ Millson and M.\ S.\ Raghunathan,
{\em Geometric construction of cohomology for arithmetic
groups I}, 
in ``Geometry and Analysis, Papers dedicated to the Memory of V.\ K.\.
Patodi, Proc. Indian Acad. Sci. vol. {\bf 90} (1981), pg. 103-123.

\bibitem[Ra]{Raghunathan}
M.\ S.\ Raghunathan,
{\em On the first cohomology of discrete subgroups of semisimple Lie groups},
Am. J. Math. vol. {\bf 87} (1965), pg. 103--139.


\bibitem[RS]{RohlfsSpeh}
J.\ Rolfs and B.\ Speh,
{\em Representations with cohomology in the discrete spectrum of subgroups
of $SO(n,1)(\mathbb{Z})$ and Lefschetz numbers},
Ann. Sci. Ec. Norm. Sup. IV vol. {\bf 20}  (1987), pg. 89--136.


\bibitem[S]{Speh}
B.\ Speh,
{\em Induced representations and the cohomology of discrete subgroups},
Duke Math. J. vol. {\bf 49} (1982), pg. 1115--1127.


\bibitem[VZ]{VoganZuckerman}
D.\ A.\ Vogan and G.\ J.\ Zuckerman,
{\em Unitary representations with non-zero cohomology},
Compositio Math. vol. {\bf 53} (1984), pg. 51--90.

\end{thebibliography}
\end{document}